\documentclass[%
]{mpi2015-cscpreprint}

\usepackage[american]{babel}

\usepackage{graphicx}

\usepackage{amssymb}
\usepackage{amsmath}
\usepackage{amsthm}
\usepackage{todonotes}
\usepackage{mathtools}
\usepackage{url}

\usepackage{algorithm}
\usepackage{algpseudocode}

\usepackage{bbm}

\usepackage{siunitx}
\usepackage{booktabs}
\usepackage{csquotes}
\usepackage{bm}

\def\IR{{\mathbb R}}
\def\IC{{\mathbb C}}

\def\IL{{\mathbb L}}

\newcommand{\bA}{{\textbf A}}
\newcommand{\bB}{{\textbf B}}
\newcommand{\bC}{{\textbf C}}
\newcommand{\bD}{{\textbf D}}
\newcommand{\bE}{{\textbf E}}

\newcommand{\bS}{{\textbf S}}
\newcommand{\bY}{{\textbf Y}}
\newcommand{\bJ}{{\textbf J}}

\newcommand{\bM}{{\textbf M}}

\newcommand{\bI}{{\textbf I}}
\newcommand{\bH}{{\textbf H}}

\newcommand{\bW}{{\textbf W}}
\newcommand{\bR}{{\textbf R}}
\newcommand{\bT}{{\textbf T}}
\newcommand{\bZ}{{\textbf Z}}
\newcommand{\bX}{{\textbf X}}
\newcommand{\bx}{{\textbf x}}
\newcommand{\by}{{\textbf y}}
\newcommand{\bu}{{\textbf u}}
\newcommand{\bV}{{\textbf V}}
\newcommand{\bU}{{\textbf U}}

\newcommand{\bfe}{{\textbf e}}

\newcommand{\bff}{{\textbf f}}

\newcommand{\bfz}{{\mathbf 0}}

\newcommand{\cC}{ {\cal C} }
\newcommand{\cD}{ {\cal D} }

\newcommand{\cN}{ {\cal N} }
\newcommand{\cH}{ {\cal H} }

\newcommand{\bPi}{ \boldsymbol{\Pi} }

\newcommand{\bomega}{ \boldsymbol{\omega} }
\newcommand{\bOmega}{ \boldsymbol{\Omega} }

\newcommand{\bDelta}{\boldsymbol{\Delta}}
\newcommand{\bLambda}{\boldsymbol{\Lambda}}

\newcommand{\bXi}{\boldsymbol{\Xi}}

\newcommand{\cL}{ {\cal L} }

\def\IR{{\mathbb R}}
\def\IC{{\mathbb C}}
\def\IL{{\mathbb L}}
\def\IV{{\mathbb V}}
\def\IW{{\mathbb W}}
\newcommand{\sIL}{{{{\mathbb L}_s}}}

\newcommand{\bone}{{\mathbbm{1}}}

\newcommand{\Sig}{\boldsymbol{\Sigma}}
\newcommand{\bom}{ \boldsymbol{\omega} }
\newcommand{\trans}{\ensuremath{^{\mkern-1.5mu\mathsf{T}}}}

\newcommand{\herm}{\ensuremath{^{\mathsf{H}}}}

\newcommand{\inv}{\ensuremath{^{-1}}}

\newcommand{\Linf}{\ensuremath{\mathcal{L}_{\infty}}}
\newcommand{\linfnorm}[1]{\ensuremath{\left\lVert#1\right\rVert_{\Linf}}}

\newcommand{\AFAAALF}{\textsf{LS-Loewner}}
\newcommand{\AAA}{\textsf{Modified AAA}}
\newcommand{\LoewnerSVD}{\textsf{Loewner-SVD}}
\newcommand{\AAACUR}{\textsf{Loewner-CUR}}
\newcommand{\lfpp}{\textsf{LFPP}}
\newcommand{\lfapp}{\textsf{LFaPP}}

\title{
Practical challenges in data-driven interpolation: dealing with noise, enforcing stability, and computing realizations
}

\author[$\ast$]{Quirin Aumann}
\affil[$\ast$]{Max Planck Instiute for Dynamics of Complex Technical Systems,
	Sandtorstr. 1, 39106 Magdeburg, Germany.\authorcr
	\email{aumann@mpi-magdeburg.mpg.de}, \orcid{0000-0001-7942-5703}}

\author[$\dagger$]{Ion Victor Gosea}
\affil[$\dagger$]{Max Planck Instiute for Dynamics of Complex Technical Systems,
	Sandtorstr. 1, 39106 Magdeburg, Germany.\authorcr
	\email{gosea@mpi-magdeburg.mpg.de}, \orcid{0000-0003-3580-4116}}

\shorttitle{Data-driven interpolation: challenges and solutions}
\shortauthor{Q. Aumann, I. V. Gosea}
\shortdate{}


\begin{document}
	
\abstract{%
	In this contribution, we propose a detailed study of interpolation-based data-driven methods that are of relevance in the model reduction and also in the systems and control communities. The data are given by samples of the transfer function of the underlying (unknown) model, i.e., we analyze frequency-response data. We also propose novel approaches that combine some of the main attributes of the established methods, for addressing particular issues. This includes placing poles and hence, enforcing stability of reduced-order models, robustness to noisy or perturbed data, and switching from different rational function representations. We mention here the classical state-space format and also various barycentric representations of the fitted rational interpolants. We show that the newly-developed approaches yield, in some cases, superior numerical results, when comparing to the established methods. The numerical results include a thorough analysis of various aspects related to approximation errors, choice of interpolation points, or placing dominant poles, which are tested on some benchmark models and data-sets.}

\novelty{This note shows that by combining the features of established data-driven rational approximation methods based on interpolation (and/or least squares fit), one can devise methods that offer additional important advantages. These include stability enforcement by placing poles in an elegant and numerically stable manner, together with robustness to noisy data.}

\maketitle

\section{Introduction}

Approximation of large-scale dynamical systems is pivotal for serving the scopes of efficient simulation and designing control laws in real-time. The technique for reducing the complexity of a system is known as model order reduction (MOR) \cite{ACA05,BOCW17,morAntBG20,benner2021system}. There exist a number of methodologies for reducing large-scale models, and each method is tailored to some specific applications (mostly, but not restricted to mechanical and electrical engineering) and to achieving certain goals (stability, passivity or structure preservation), on top of the complexity reduction part. Data-driven MOR approaches are of particular importance when access to high-fidelity models is not explicitly granted. This means that a state-space formulation with access to internal variables is not available, yet input/output data are.
Such methods circumvent the need to access an exact description of the original model and are applicable whenever classical projection-based MOR is not. Here, we mention system and control methodologies that are based on interpolation or least-square fit of data (i.e., frequency response measurements), such as vector fitting \cite{morGusS99}, the Loewner framework \cite{MA07}, or the AAA algorithm \cite{NST18}. Methods that use time-domain data are also of interest, including the ones that require input-output data together with those which use snapshot data (access to the state evolution), such as the classical ones in \cite{JuangPappa85,van1994n4sid,van2012subspace}, followed by \cite{morAst10}, \cite{PGW17} or \cite{schmid_2010}.

We focus on interpolation-based or so-called moment matching (MM) methods which have emerged, were developed, and improved continuously in the last decades. The backbone of such methods represent rational Krylov-type approaches together with the Sylvester matrix equation interpretation \cite{bai2002krylov,gallivan2004sylvester,ACA05}. Apart from being computationally efficient and easy to implement, MM approaches have another advantage: they do not (necessarily) require access to a full-state realization of the original dynamical system. Hence, they can be viewed as data-driven methods. Here data are given by the moments of the system, i.e., samples of the underlying transfer function of the system (and of its derivatives) evaluated in a particular frequency range; for more details, we refer the readers to \cite{MA07,morAst10,morAntBG20} and to Chapter~3 in \cite{benner2021system}. The notion of a moment with respect to systems and control theory is related to the unique solution of a Sylvester matrix equation \cite{gallivan2004sylvester}.

The purpose of this note is twofold; first, we intend to review and to connect three important system theoretical model reduction methods based on interpolation that were introduced in the last 15 years:
\begin{itemize}
	\item  The Loewner framework (LF) by Mayo and Antoulas from 2007 in \cite{MA07};
	\item The Astolfi framework (AF) from 2010 in \cite{morAst10};
	\item  The Adaptive Antoulas Anderson (AAA) algorithm by Nakatsukasa, Set\'e and Trefethen from 2018 in \cite{NST18}.
\end{itemize}
Together, these three approaches were cited multiple times in various research publications, being arguably quite popular methods. However, until now, not too many connections between them were provided, neither in the automatic control, nor in the model reduction, or numerical analysis communities. Together with the vector fitting algorithm (VF) in \cite{morGusS99} (which is not based on interpolation, and is hence a purely optimization approach based on least-squares fitting), these methods represent arguably the most prolific rational approximation schemes developed in the system and control community.  However, VF is not the object of this study since it is not based on interpolation.

The other scope of this note is to propose a new method that is based on the three methods enumerated above, and that addresses some of the shortcomings and challenges associated with them.
Basically, the idea is to combine the attributes of each method, by following the steps below.
\begin{itemize}
	\item We make use of the order-revealing property of the LF (encoded by the rank of augmented Loewner matrices); additionally, the selection of interpolation points is done via a Loewner-CUR technique proposed in \cite{morKarGA21a}.
	\item We utilize the elegant state-space parameterization of the LTI system proposed by the AF (after imposing $k$ interpolation conditions); this is the backbone of the methods (we also show the connection between state-space forms and barycentric forms).
	\item We use either the fitting step from AAA (that chooses free parameters to fit the un-interpolated data in a least square sense) or we impose pole placing (dominant poles are selected from those of the Loewner model); in both cases, a linear system of equations needs to be solved. 
\end{itemize}

\noindent
In what follows, we consider a multiple-input multiple-output (MIMO) linear time-invariant (LTI) system $\Sig_L$ of dimension $n$ described by the following system of differential equations:
\begin{equation}\label{def_lin_sys}
	\Sig_L: \begin{cases}
		\dot{\bx}(t) = \bA \bx(t) + \bB \bu(t), \ \
		\by(t)  = \bC \bx(t),
	\end{cases}
\end{equation} 
with $\bx(t) \in \IR^n$ as the state variable, $\bu(t) \in \IR^m$ as the control inputs, and $\by(t) \in \IR^p$ as the observed outputs. Here, we have that $\bA \in \IC^{n \times n}$, $\bB \in \IC^{n \times m}$ and $\bC \in \IC^{p \times n}$. The transfer (matrix) function of the LTI system is given by $\bH(s) \in \IC^{p \times n}$, with $s \in \IC$, as
\begin{equation}\label{eqn:Htrf}
	\bH(s) = \bC (s \bI_n - \bA)^{-1} \bB.
\end{equation}
It is to be noted that, for $m = p =1$, the system becomes  single-input single-output (SISO). We will sometime switch between MIMO and SISO formats while presenting the methods covered in this note, since the latter allows a more easy exposition for some of the results shown here.

Let $s_i \in \mathbb{C} \setminus \sigma(\bA)$, where $\sigma(\bA)$ denotes the spectrum of matrix $\bA  \in \IC^{n \times n}$, i.e., the set of its eigenvalues. The $j$-moment of system $\Sig_L$ at $s_i$ is given by $\eta_j(s_i) = \frac{(-1)^j}{j!} \Big{[} \frac{d^j}{ds^j} \bH(s) \Big{]}_{s=s_i}$, for any integer $j \geqslant 1$. The 0-moment is obtained by sampling the transfer function $\bH(s)$ in \cref{eqn:Htrf} at $s_i$, i.e., $\eta_0 = \bH(s_i)$. 
In this contribution, we restrict the analysis to matching 0-moments, i.e., samples of the transfer function $\bH(s)$, and not of its derivatives. However, all methodologies shown here can be expected to cope with this as well. Moreover, in practice, inferring 0-moments from time-domain data is usually a more straight-forward task; this is performed by exciting the system with harmonic inputs, and by applying spectral transformations to the outputs. Additionally, the inference of derivative values (of the transfer function) is typically susceptible to perturbations and it is more challenging to attain, from a numerical point of view.

The paper is structured in the following way; after the introduction session sets up the stage, we propose a survey of three established interpolation-based methods in \Cref{sec:methods}. Then, the proposed methodologies are developed in \Cref{sec:main}, with emphasis on the one-step approach that combines optimal selection of interpolation points (chosen using CUR-DEIM) with LS fit on the rest of the data, and also the pole placement method in barycentric form that enforces dominant poles from the Loewner data-driven model. Then, \Cref{sec:numerics} illustrates the numerical aspects of applying the methods  discussed/proposed in the previous two sections to a variety of test cases (various models and data sets). Finally, \Cref{sec:conc} presents the conclusions and the outlook into future research.

\section{A survey of established methods}
\label{sec:methods}

In this section we discuss three established data-driven methods for rational approximation (AF, LF, and AAA, as mentioned in the previous section). The data are samples of the transfer function corresponding to the underlying dynamical system, measured on a particular frequency grid. In what follows, we mention some state-of-the-art methodologies used to measure such data, i.e., frequency response data. Typically, such measurements can be produced in practice from experiments conducted in scientific laboratories using carefully calibrated machines, called spectrum analyzers (SAs). In this category we mention swept-tuned spectrum analyzers, scalar network analyzers (SNAs), and vector network analyzers (VNAs). 

The SNA is an instrument that measures microwave
signals by converting them to a DC voltage using a diode detector.
In a VNA, information regarding both the magnitude and the phase of a microwave signal is extracted. While there are different ways to perform such measurements, the method employed by commercial products (such as the Anritsu series described in \cite{Anritsu}) of VNAs is to down-convert the signal to a lower intermediate frequency in a process called harmonic sampling. This signal can then be measured directly by a tuned receiver. Compared to the SNA, the VNA is a more powerful analyzer tool. The major difference is that the VNA can also measure the phase, and not only the amplitude. With this property enforced, then so-called scattering parameters (or S-parameters) can be processed. These can be used for identifying forward and reverse transmission and reflection characteristics. More details can be found in \cite{Anritsu}.

The harmonic balance method (or HBM) \cite{NV76}, is an established methodology in the field of electromagnetics. The HBM is used in many (if not most) commercial radio-frequency (RF) simulation tools. This is due to the fact that it has certain advantages over other common methods, namely modified nodal analysis (MNA), which makes it more appropriate to use for stiff problems and circuits containing
transmission lines, nonlinearities and dispersive effects.  More details can be found in the survey paper \cite{peyton2018recent}.

\subsection{The one-sided moment-matching approach in \cite{morAst10}}

The framework introduced by Astolfi in \cite{morAst10} (referred to as AF throughout the paper) deals with the problem of
model reduction by moment matching. Although classically
interpreted as a problem of interpolation of points in the
complex plane, it has instead been recast as a problem of interpolation of steady-state responses. In the following we briefly review its application to linear systems. It is to be noted that the AF was steadily extended and applied to different scenarios (including nonlinear dynamical systems, pole-zero placement, and least-squares fit) \cite{i-astolfi-colaneri-SCL2014,i-astolfi-TAC2016,morScaA17,simard-astolfi-TAC2021,padoan2021model,Ionescu2022}.

The moments of a linear system can be characterized in terms of the solution of Sylvester equations. By using this observation, it has been shown that the moments are in one-to-one relation with the steady-state output response of the interconnection between a signal generator and the  original linear system.

In what follows, for simplicity of exposition, it is assumed that $\Sig_\mathrm{L}$ is a minimal system (both fully controllable and fully observable). For exact definitions on minimality, controllability, or observability of LTI systems, we refer the reader to \cite{ACA05}.

Let $k \leqslant n$ and $\bS \in \mathbb{C}^{k \times k}$ a non-derogatory matrix (for which the characteristic and minimal polynomials coincide) with $\sigma(\bS) \cap \sigma(\bA) = \emptyset$ and $\bR \in \mathbb{C}^{1 \times k}$ so that $(\bS,\bR)$ is observable. Consider the signal generator system $\Sig_{\textrm{sg}}$ described by the equations
\begin{equation}\label{def_sig_gen}
	\Sig_{\textrm{sg}}: \begin{cases}
		\dot{\bom}(t)  = \bS \bom(t) , \ \
		\bu(t)  = \bR \bom(t).
	\end{cases}
\end{equation}
Then, the explicit solution of (\ref{def_sig_gen}) can be written as $\bom(t) = e^{\bS t} \bom(0)$. Hence, the control input is written as $\bu(t) = \bR e^{\bS t} \bom(0)$. In addition, the eigenvalues of $\bS$ are called interpolation points.

For a linear system $\Sig_L$, and interpolation points $s_i \in \mathbb{C} \setminus \sigma(\bA)$, for $i = 1,\ldots,k$, consider a non-derogatory matrix $\bS \in \mathbb{R}^{k \times k}$. It follows that there exists a one-to-one relation between the moments of the system $\Sig_L$ and 
\begin{enumerate}
	\item the matrix $\bC \bPi$, where $\bPi$ is the (unique) solution of the Sylvester equation $\bA \bPi +\bB \bR = \bPi \bS$, for any row vector $\bR \in \mathbb{R}^{1 \times k}$ so that $(\bR,\bS)$ is observable;
	\item the steady-state response of the output $\by$ of the interconnection of system $\Sig_L$ and the system $\Sig_\mathrm{sg}$, for any $\bR$ and $\bom(0)$ such that the triplet $\left(\bR,\bS,\bom(0)\right)$ is minimal.
\end{enumerate}
More precisely, let $\bDelta \in \mathbb{R}^k$ be a column vector containing $k$ free parameters (denoted here by $\delta_1, \delta_2, \ldots, \delta_k$, with $\delta_i \neq 0, \ 1 \leq i \leq k$). Then, as stated in \cite{morAst10}, the family of linear time-invariant systems that interpolates the moments of system $\Sig_{\textrm{L}}$ at the eigenvalues of matrix $\bS$, is given by
\begin{equation}\label{Astolfi_right_model}
	\hat{\Sig}_{\bDelta}: \begin{cases}
		\dot{\hat{\bx}}(t) = \underbrace{(\bS-\bDelta \bR)}_{=\hat{\bA}} \hat{\bx}(t) + \underbrace{\bDelta}_{=\hat{\bB}} \bu(t), \ \ \hat{\by}(t) = \underbrace{\bC \bPi}_{=\hat{\bC}} \hat{\bx}(t),
	\end{cases}
\end{equation} 
where the matrices $\bS$ and $\bR$ are as before and $\bPi$ is the unique solution of the Sylvester equation $\bA \bPi +\bB \bR = \bPi \bS$. Additionally, the condition $\sigma(\bS) \cap \sigma(\bS-\bDelta \bR) = \emptyset $ needs to be enforced. It is to be noted that the free parameters explicitly enter the vector $\hat{\bB} = \bDelta$, but also the matrix $\hat{\bA}$, as $\hat{\bA} = \bS-\bDelta \bR$. Finally, $\hat{\bC} = \bC \bPi$ has fixed entries.

The Sylvester matrix equation for the reduced-order system is written as $\hat{\bA} \hat{\bPi} +\hat{\bB} \bR = \hat{\bPi} \bS$.  This can be explained by the fact that the reduced-order system matches the prescribed moments of the original system, hence the same format of the two equations. Without loss of generality, one can consider that the matrix $\hat{\bPi}$ is the identity matrix, i.e. $\hat{\bPi} = \bI_k$ (this can be achieved by applying similarity transformations). By replacing this value into the reduced-dimension Sylvester matrix equation above, the formula $\hat{\bA} = \bS-\bDelta \bR$ directly follows.

Afterwards, the free parameters collected in the vector $\bDelta$ can be chosen in order to enforce or impose additional conditions as mentioned in \cite{morAst10}, such as: matching with imposing additional $k$ interpolation conditions, matching with prescribed eigenvalues, matching with prescribed relative degree, matching with prescribed zeros, matching with a passivity constraint, matching with $L_2$-gain, or matching with a compartmental constraint.

An important aspect of the AF is the characterization of all, i.e., infinitely many families of reduced-order models that satisfy $k$ prescribed interpolation conditions. This is done by explicitly computing such parameterized models, for which the free parameters are the variables entering in the vector $\bDelta$. The main parameterization developed here will be used as a \enquote{backbone} of the methods developed in \Cref{sec:main}.

As stated in the original paper, the main advantage of the AF (characterization of moments in terms of steady-state responses) is that it allows the definition of moments for systems which do not admit a clear/immediate representation in terms of transfer function(s). Hence, the author provides as examples the case of linear time-varying systems, and the case of nonlinear systems. Moreover, it is stated in \cite{morAst10} that one disadvantage of the framework is that it requires the existence of steady-state responses. Consequently, the original system has to be (exponentially) stable. However, in most practical applications, this is a realistic requirement.

\subsection{The Loewner framework in \cite{MA07}}
\label{sec:LF}

In this section we present a short summary of the Loewner framework (LF), as introduced in \cite{MA07}. It is to mentioned that LF has its roots in the earlier work of \cite{AA86}, and that LF can be considered to be a double-sided moment-matching approach (as opposed to AF, which is one-sided). 

For a tutorial paper on LF for LTI systems, we refer the reader to \cite{ALI17}, and for a recent extension that uses time-domain data, we refer the reader to \cite{morPehW16}. 
The Loewner framework has been recently extended to certain classes of nonlinear systems, such as bilinear systems in \cite{AGI16}, and quadratic-bilinear (QB) systems in \cite{GA18}, but also to linear parameter-varying systems in \cite{GPA21CDC}. Additionally, issues such as stability preservation or enforcement, or passivity preservation, were tackled in the LF in \cite{morGosA16,morGosPA21a}, for the former, and in \cite{antoulas2005new, benner2020identification}, for the latter.

The LF is based on processing frequency-domain measurements
$\cD = \{\left(\omega_\ell,\bH(\omega_\ell)\right), ~\ell=1,\ldots,N\}$ (with $\omega_\ell \in \IR$ for $1 \leq \ell \leq N$) corresponding to evaluations of the transfer function of the underlying (unknown/hidden) dynamical system. 

The interpolation problem is formulated as shown below (for convenience of exposition, we show here only the SISO formulation). We are given data nodes and data points in the set $\cD$, partitioned into two disjoint subsets $\cD_{\textrm L}$ and $\cD_{\textrm R}$, with $\cD_{\textrm L} \cup \cD_{\textrm R} = \cD$ and $k+q = N$, as
\begin{align}\label{data_Loew}
	\begin{split}
		{\textrm{right \ data}}&: \cD_{\textrm L} = \{\left(\lambda_j,\bH(\lambda_j)\right), ~j=1,\ldots,k\},~{\textrm {and}}, \\
		{\textrm{left \ data}}&: \cD_{\textrm R}=\{\left(\mu_i,\bH(\mu_i)\right), ~i=1,\ldots,q\},
	\end{split}
\end{align}
and we seek to find a rational function
$\hat{\bH}(s)$, such that the following interpolation conditions hold:
\begin{equation} \label{interp_cond}
	\hat{\bH}(\mu_i)=\bH(\mu_i):=v_i,~~~\hat{\bH}(\lambda_j)=\bH(\lambda_j):= w_j.
\end{equation}
The Loewner matrix $\IL \in\IC^{q\times k}$ and the shifted Loewner matrix $\sIL \in\IC^{q\times k}$ play an important role in the LF, and are given by
\begin{equation} \label{Loew_mat}
	\IL_{(i,j)}=\frac{v_i-w_j}{\mu_i-\lambda_j}, \ \ \sIL_{(i,j)}=
	\frac{\mu_i v_i-\lambda_j w_j}{\mu_i-\lambda_j},
\end{equation}
while the data vectors $\IV \in \IC^q,\ \IW\trans \in \IC^k$ are given by
\begin{equation} \label{VW_vec}
	\IV_{(i)}= v_i, \ \  \IW_{(j)} = w_j,~\text{for}~i=1,\ldots,q, \ j=1,\ldots,k.
\end{equation}
Moreover, the following Sylvester matrix equations (\cite[Ch. 6]{ACA05}) are satisfied by the Loewner and shifted Loewner matrices (here, $\bone_q = \left[ \begin{matrix} 1 &  \cdots & 1  \end{matrix} \right]\trans \in \IC^q$)
\begin{equation}\label{eq:sylv_loe}
\begin{cases} \bM \IL - \IL \bLambda = \IV \bone_k\trans - \bone_q \IW, \\
\bM \sIL - \sIL \bLambda =  \bM \IV \bone_k\trans -  \bone_q  \IW \bLambda,
\end{cases}
\end{equation}
where $\bM = \text{diag}(\mu_1,\ldots,\mu_q)$ and $\bLambda = \text{diag}(\lambda_1,\ldots,\lambda_k)$ are diagonal matrices. The following relation holds true
\begin{equation}
\sIL =  \IL \bLambda + \IV \bone_k\trans = \bM \IL + \bone_q  \IW.
\end{equation}
The unprocessed Loewner surrogate model, provided that $k=q$, is composed of the matrices
\begin{align}\label{Loew_doub_sided}
	\hat{\bE}=-\IL,~~ \hat{\bA}=-\sIL,~~ \hat{\bB}=\IV,~~ \hat{\bC}=\IW,
\end{align}
and if the pencil $(\IL,\sIL)$ is regular, then the function $\hat{\bH}(s)$ satisfying the interpolation conditions in \cref{interp_cond} can be explicitly computed in terms of the matrices in \cref{Loew_doub_sided}, as $\hat{\bH}(s) = 	\hat{\bC} (s	\hat{\bE} - \hat{\bA})^{-1} \hat{\bB}$.

In practical applications (when processing a fairly large number of measurements), the pencil $(\sIL,\,\IL)$ is often singular. Hence, a post-processing step is required for the Loewner model in \cref{Loew_doub_sided}. In such cases, one needs to perform a singular value decomposition (SVD) of augmented Loewner matrices, to extract the dominant features and remove inherent redundancies in the data. By doing so, projection matrices $\bX, \bY \in \IC^{k \times r}$ are obtained, as left, and respectively, right truncated singular vector matrices:
\begin{equation}
	\left[\IL~~\IL_s\right]=\bY\bS_{ {r}}^{(1)}\tilde{\bX}\herm~\left[\begin{array}{c}\IL \\ \IL_s\end{array}\right] = {\tilde\bY}\bS_{ {r}}^{(2)} \bX\herm,
\end{equation}	
where $\bS^{(1)}_{ {r}}$, $\bS_{ {r}}^{(2)}$ $\in\IR^{{{r}}\times{r}}$,~
$\bY \in\IC^{k\times{r}}$,\ $\bX\in\IC^{q\times{r}}$,~$\tilde{\bY}\in\IC^{2q\times{r}}$,~$\tilde{\bX}\in\IC^{r\times{2k}}$. The truncation index $r$ can be chosen as the \textit{numerical rank} (based on a tolerance value $\tau >0$) or as the \textit{exact rank} of the Loewner pencil (in exact arithmetic), depending on the application and data size. More details can be found in \cite{ALI17}.

The system matrices corresponding to a projected Loewner model of dimension $r$ can be computed as follows:
\begin{equation*}
	\tilde{\bE} = -\bX\herm\IL \bY, \ \ \tilde{\bA} = -\bX\herm\sIL \bY, \ \
	\tilde{\bB} = \bX\herm\IV, \ \  \tilde{\bC} = \IW \bY.
\end{equation*}
We note that MIMO extensions of the LF were already proposed in the original contribution \cite{MA07}. There, a tangential interpolation framework is considered. Instead of imposing interpolation of full $p \times m$ blocks, the authors prefer to interpolate the original transfer matrix function samples along certain vectors (or tangential directions). We also note that a first attempt of re-interpreting the LF in \cite{MA07} as a one-sided method was  made in \cite{gosea2022one}. In the latter, the main difference to the classical work in \cite{AA86} was that a compression of the left (un-interpolated) data set was enforced. However, in \cite{gosea2022one}, it was still unclear how to split the data, i.e., what the right data set should be (where interpolation is enforced). Finally, it is to be noted that the choice of interpolation points is crucial in the LF. An exhaustive study of different choices was proposed in \cite{morKarGA19a}, while a greedy strategy was proposed in \cite{cherifi2022greedy}, for scenarios in which limited experimental data are available.

\subsection{The AAA algorithm in \cite{NST18}}
\label{sec:AAA}

The AAA algorithm introduced in \cite{NST18} is an adaptive and iterative extension of the interpolation method based on Loewner matrices, originally proposed in \cite{AA86}. The main steps are as follows
\begin{enumerate}
	\item Express the fitted rational approximants in a barycentric representation, which represents a numerically stable way of expressing rational functions \cite{berrut2004barycentric}.
	
	\item Select the next interpolation (support) points via a greedy scheme; basically, interpolation is enforced at the point where the (absolute or relative) error at the previous step is maximal.

	\item Compute the other variables (the so-called barycentric weights) in order to enforce least squares approximation on the non-interpolated data.
\end{enumerate}	

In recent years, the AAA algorithm has proven to be an accurate, fast, and reliable rational approximation tool with a fairly large range of applications. Here, we will mention only a few: nonlinear eigenvalue problems \cite{lietaert2022automatic}, MOR of parameterized linear dynamical systems \cite{CS20}, MOR of linear systems with quadratic outputs \cite{gosea2022data}, rational approximation of periodic functions \cite{baddoo2021aaatrig}, representation of conformal maps \cite{gopal2019representation}, rational approximation of matrix-valued functions \cite{GG20}, or signal processing with trigonometric rational functions \cite{wilber2022data}. The procedure is sketched in \Cref{al:aaa}.

\begin{algorithm}[tb]
	\caption{The AAA algorithm.}
	\label{al:aaa}
	\begin{algorithmic}[1]
		\Require A (discrete) set of sample points $\Gamma \subset \mathbb{C}$ with $N$ points, function $f$ (or the evaluations of $f$ on the set $\Gamma$, i.e., the sample values), and an error tolerance $\epsilon>0$.
		\Ensure A rational approximant $r_{n}(s)$ of order $(n,n)$ displayed in a barycentric form.
		\State Initialize $j=0$, $\Gamma^{(0)} \gets \Gamma$, and $r_{-1} \gets N^{-1}\sum_{i=1}^N f(\gamma_i)$.
		\While{$\vert f(s)-r_{j-1}(s) \vert > \epsilon$}
		\State \parbox[t]{\dimexpr\textwidth-\leftmargin-\labelsep-\labelwidth}{%
			Select a  point $z_{j} \in \Gamma^{(j)}$ for which $\vert f(s)-r_{j-1}(s) \vert$ attains a maximal value, where for $j \geq 1$, it follows:}
		\begin{equation}\label{eqn:proper_trsfct}
			r_{j-1}(s) :=  \left(\displaystyle \sum_{k=0}^{j-1}\frac{\omega_k^{(j-1)}}{s - z_k}\right)^{-1} \left(\displaystyle \sum_{k=0}^{j-1} \frac{\omega_k^{(j-1)} f_k}{s - z_k}\right).
		\end{equation}
		\If{$\vert f(z_j) - r_{j-1}(z_j)\vert \leq \varepsilon$}
		\State Return $r_{j-1}$.
		\Else 
		\State $f_j \gets f(z_j)$ and $\Gamma^{(j+1)} \gets \Gamma^{(j)}\setminus \{z_j\}$.
		\EndIf
		\State Find the weights $\bomega^{(j)} = [\omega_0^{(j)},\ldots,\omega_j^{(j)}]$ by solving a least squares problem over $z\in\Gamma^{(j+1)}$
		\begin{align}\label{eqn:LSfit_Loew}
			\begin{split}
				&\displaystyle \sum_{k=0}^j\frac{\omega_k^{(j)}}{s - z_k} f(s) \approx \sum_{k=0}^j \frac{\omega_k^{(j)} f_k}{s - z_k} \Leftrightarrow \ \  \Big{(} \sum_{k=0}^j\frac{f(s) - f_k}{s - z_k} \Big{)} \omega_k^{(j)} \approx 0  \Leftrightarrow \IL^{(j)} \bomega^{(j)} = 0.
			\end{split}
		\end{align}
		\hskip\algorithmicindent \parbox[b]{\dimexpr\textwidth-\leftmargin-\labelsep-\labelwidth}{%
			The solution of \cref{eqn:LSfit_Loew} is given by the $(j+1)^\mathrm{th}$ right singular vector of the Loewner matrix $\IL^{(j)} \in \mathbb{C}^{(N-j-1)\times(j+1)}$.}
		\State $j \gets j+1$.
		\EndWhile
	\end{algorithmic}
\end{algorithm}

It is to be mentioned that a modified version of AAA that enforces \textbf{real-valued} and  \textbf{strictly-proper} rational appoximants was recently proposed in \cite{gosea2021loewner}. There, the format of the function in \cref{eqn:proper_trsfct} was modified by inserting a $1$ into the denominator, as follows
	\begin{equation}\label{eqn:str_proper_trsfct}
\tilde{r}_{j}(s) :=  \left(\displaystyle 1+\sum_{k=0}^{j-1}\frac{\omega_k^{(j-1)}}{s - z_k}\right)^{-1} \left(\displaystyle \sum_{k=0}^{j-1} \frac{\omega_k^{(j-1)} f_k}{s - z_k}\right).
\end{equation}
Consequently, the equation in \cref{eqn:LSfit_Loew} becomes $\IL^{(j)} \bomega^{(j-1)} = -\bff^{(j-1)}$, where the vector $\bff^{(j-1)} \in  \IC^j$ is given by $\bff^{(j-1)} = \left[ \begin{matrix} f_0 & f_2 & \cdots & f_{j-1} \end{matrix} \right]\trans$. It is to be noted that $\tilde{r}_{j}(s)$ in \cref{eqn:str_proper_trsfct} is theoretically a rational approximant of order $(j-1,j)$, if we do not take into account pole/zero cancellations or any other zero cancellations of coefficients in the numerator or denominator.
\section{The proposed methodologies}
\label{sec:main}
\subsection{Skeleton of the main methods}
\label{sec:skeleton}
Similar to the methods reviewed in \Cref{sec:methods} we want to find an LTI system with a transfer function of the structure \cref{def_lin_sys} that interpolates data provided as measurements $\bH\left(s_i\right),\ i=1, \ldots, k$ of the transfer function of the original system.
We can directly put together an LTI parametrized model of dimension $r=km$, having $km^2$ degrees of freedom
with transfer function
\begin{equation}\label{eqn:Htrf_r}
	\hat{\bH}(s) = \hat{\bC} (s \bI_r - \hat{\bA})^{-1} \hat{\bB},
\end{equation}
with the underlying data concatenated to
\begin{equation}
	\label{eq:realization_C}
	\hat{\bC} = \begin{bmatrix} \bH(\lambda_1) & \cdots & \bH(\lambda_k)
\end{bmatrix} \in \IC^{p \times r},
\end{equation}
a matrix of weights $\hat{\bW}_i$
\begin{equation}
	\label{eq:realization_B}
	\hat{\bB} = \begin{bmatrix} \hat{\bW}_1\herm & \cdots & \hat{\bW}_k\herm
	\end{bmatrix}\herm \in \IC^{r \times m},
\end{equation}
and $\hat{\bA} \in \IC^{r \times r}$ formed from a diagonal matrix populated with the interpolation points $\lambda_i$ disturbed by $\hat{\bB}$, such that
\begin{equation}
	\label{eq:realization_A}
	\hat{\bA} = \bLambda - \hat{\bB} \bR = \operatorname{diag}\left( \lambda_1, \ldots, \lambda_k\right)
	\otimes \bI_m - \hat{\bB} \left(\bone_r\trans \otimes \bI_m\right).
\end{equation}
Making use of the Woodbury matrix identity
and denoting $\bLambda_s  = s\bI_{km}-\bLambda$, the transfer function \cref{eqn:Htrf_r} can be rewritten as
\begin{equation}
	\label{eq:rom_woodbury}
	\hat{\bH}(s) = \hat{\bC} \bLambda_s^{-1} \hat{\bB} \left( \bI_m + \bR \bLambda_s^{-1} \hat{\bB} \right)^{-1}.
\end{equation}
A complete derivation of \cref{eq:rom_woodbury} is given in \Cref{sec:app_wmi}.

In the single-input single-output case ($m = p =1$, hence $r=k$), the barycentric weights reduce to scalars and the matrices for a ROM of structure \cref{eqn:Htrf_r} are given by
\begin{align}\label{ROM_SISO}
\begin{split}
\hat{\bA} &= \bLambda - \hat{\bB} \bR  \in \IC^{k \times k}, \ \ \hat{\bB} = \begin{bmatrix} \hat{w}_1 & \cdots & \hat{w}_k \end{bmatrix}\trans \in \IC^{k \times 1}, \ \ \hat{\bC} = \begin{bmatrix} \bH\left(\lambda_1\right)& \cdots & \bH\left(\lambda_k\right)
\end{bmatrix} \in \IC^{1 \times k}.
\end{split}
\end{align}
By inserting the formulae in \cref{ROM_SISO} into \cref{eq:rom_woodbury}, and using the notation $h_i := \bH\left(\lambda_i\right)$, leads to
\begin{align}
\hat{\bC} \bLambda_s^{-1} \hat{\bB} &= \sum_{i=1}^k \displaystyle \frac{\hat{w}_i h_i}{s-\lambda_i},\ \
\left(\bI_m+ \bR \bLambda_s^{-1} \hat{\bB} \right)^{-1} =  \frac{1}{1+\sum_{i=1}^k \displaystyle \frac{\hat{w}_i}{s-\lambda_i}}.
\end{align}
Hence, the transfer function of the model in \cref{ROM_SISO} is given in barycentric representation by
\begin{align}
\hat{\bH}(s) &=  \frac{\sum_{i=1}^k \displaystyle \frac{\hat{w}_i h_i}{s-\lambda_i}}{1+\sum_{i=1}^k \displaystyle \frac{\hat{w}_i}{s-\lambda_i}}.
\end{align}
This can be performed analogously for a multi-input multi-output case ($m=p$, $r=km$).
The first part of \cref{eq:rom_woodbury} becomes
\begin{equation}
\hat{\bC} \bLambda_s^{-1} \hat{\bB} =
\begin{bmatrix} \bH(\lambda_1) \bI_m (s - \lambda_1)^{-1} & \cdots & \bH(\lambda_k) \bI_m (s - \lambda_k)^{-1} \end{bmatrix}
\begin{bmatrix} \hat{\bW}_1 \\ \vdots \\ \hat{\bW}_k \end{bmatrix}
= \sum_{i=1}^k \displaystyle \frac{\bH(\lambda_i) \hat{\bW}_i}{s-\lambda_i},
\end{equation}
the second part
\begin{gather}
\begin{aligned}
\left( \bI_m + \bR \bLambda_s^{-1} \hat{\bB} \right)^{-1} &=
\left( \bI_m +\begin{bmatrix} \bI_m & \cdots & \bI_m \end{bmatrix}
\begin{bmatrix} \bI_m (s - \lambda_1)^{-1}  &  &  \\  & \ddots & \\ & & \bI_m (s - \lambda_k)^{-1} \end{bmatrix}^{-1}
\begin{bmatrix} \hat{\bW}_1 \\ \vdots \\ \hat{\bW}_k \end{bmatrix}  \right)^{-1}  \\
&= \left( \bI_m + \sum_{i=1}^k \displaystyle \frac{\hat{\bW}_i  }{s-\lambda_i} \right)^{-1}.
\end{aligned}
\end{gather}
Consequently, the transfer function in \cref{eq:rom_woodbury} has also a barycentric form given by
\begin{align}\label{eq:MIMO_Bary_StrProp}
\hat{\bH}(s) &= \left( \sum_{i=1}^k \displaystyle \frac{\bH(\lambda_i) \hat{\bW}_i}{s-\lambda_i} \right) \left( \bI_m + \sum_{i=1}^k \displaystyle \frac{\hat{\bW}_i }{s-\lambda_i} \right)^{-1}.
\end{align}

The transfer function is defined by the choice of the interpolation points and of the weights. The interpolation points can be chosen as dominant parts of the available data or based on their location in the frequency spectrum. The weights can be computed such that the data which are not interpolated, are approximated in an optimal way. Alternatively, the weights can be chosen to enforce poles at specific locations. In the following, we show different strategies for both choices.

\subsection{Automatic choice of interpolation points}
The approximation quality of a surrogate model of the form \cref{eqn:Htrf_r} is greatly influenced by the choice of the interpolation points $\lambda$. This choice is not always obvious, so automatic strategies are frequently employed. The Loewner framework uses the SVD to identify dominant subsets of the available data to enforce interpolation on. Alternatively, the AAA algorithm uses a greedy scheme to minimize the error between surrogate and original data. Another approach, originally introduced by \cite{morKarGA19a}, makes use of the CUR decomposition to extract interpolation points from a relevant subset of the available data.

The CUR decomposition approximates a matrix $\bA$ by a product of three low-rank matrices $\check{\bA} = \check{\bC} \check{\bU} \check{\bR}$, where $\check{\bC}$ and $\check{\bR}$ represent subsets of the columns respectively rows of $\bA$ \cite{mahoney2009cur,sorensen2016deim}.
In our case the three matrices are only a byproduct, we are more interested in the interpolation points $\lambda$ and $\mu$ that are associated to the columns and rows extracted as $ \check{\bC}$ and $ \check{\bR}$.
In combination with the skeleton for a realization described in \Cref{sec:skeleton}, \Cref{al:ls_loewner} computes a surrogate model approximating a set of given transfer function data. We use the algorithm from \cite{sorensen2016deim} to compute the CUR decomposition and thus identify dominant parts of the original data set and their corresponding left and right interpolation points.
Contrary to \cite{morKarGA19a}, we decompose the original Loewner matrix $\IL$ rather than the augmented Loewner matrices $\begin{bmatrix}	\IL & \IL_s \end{bmatrix} $ and $\begin{bmatrix}	\IL\herm & \IL_s\herm \end{bmatrix}\herm $.
Using all interpolation points obtained from the CUR decomposition would introduce redundant data into the surrogate. Therefore we choose only a subset of the interpolation points: either only the left points, only the right points, or every other entry from a concatenated and sorted vector of left and right points.
Together with the data associated to the chosen interpolation points they are used to populate a rectangular Loewner matrix. We now need to compute weights for barycentric interpolation as described in the following section. After having obtained the weights, a surrogate model \cref{eqn:Htrf_r} can be computed from \cref{eq:realization_A,eq:realization_B,eq:realization_C}.

\begin{algorithm}[tb]
\caption{\AFAAALF{} with CUR.}
\label{al:ls_loewner}
\begin{algorithmic}[1]
	\Require Transfer function samples $\left\{\bH\left(s_i\right)\right\}_{i=1}^N$, corresponding sampling points $\bXi = \left\{s_i\right\}_{i=1}^N$.
	\Ensure Surrogate model $\hat{\bH}(s) = \hat{\bC} (s \bI_r - \hat{\bA})^{-1} \hat{\bB}$.
	\State Partition data and compute Loewner matrix $\IL$ as in \cref{Loew_mat}.
	\State Compute CUR decomposition, such that $\IL = \check{\bC} \check{\bU} \check{\bR}$ with $\check{\bC} \in \IC^{N \times k}, \check{\bR} \in \IC^{k \times N}$.
	\State Obtain interpolation points $\left\{\lambda_i\right\}_{i=1}^k, \left\{\mu_i\right\}_{i=1}^k$ corresponding to the columns and rows in $\check{\bC}, \check{\bR}$.
	\State Postprocess interpolation points to obtain $\bm{\nu} = \left\{\nu\right\}_{i=1}^k$ and $\bm{\chi} = \bXi \setminus \bm{\nu}$.
	\State Populate a rectangular Loewner matrix $\IL_{\left(i,j\right)} = \frac{\bH\left(\chi_i\right) - \bH \left(\nu_j\right)}{\chi_i - \nu_j} $.
	\State Compute the weights $\bOmega = -\IL^{\dagger} \begin{bmatrix}
	\bH\left(\nu_1\right)\herm & \cdots & \bH\left(\nu_k\right)\herm
	\end{bmatrix}\herm$, where $\IL^{\dagger}$ is the pseudo-inverse of $\IL$ and $\bOmega =  \begin{bmatrix}
	\hat{\bW}_1\herm & \cdots & \hat{\bW}_k\herm
	\end{bmatrix}\herm$.
	\State Compute $\hat{\bA}, \hat{\bB}, \hat{\bC}$ with \cref{eq:realization_A,eq:realization_B,eq:realization_C}.
\end{algorithmic}
\end{algorithm}

\subsection{Computing the barycentric weights}
\subsubsection{Least-squares approach}

The matrix-valued weights $\hat{\bW}_i$ can be computed similarly to AAA \cite{GG20} by solving the minimization problem
\begin{align}
\min_{\hat{\bW}_i} \sum_{j=1}^h \left( \left( \sum_{i=1}^k \displaystyle \frac{\bH(\lambda_i) \hat{\bW}_i}{s_j-\lambda_i} \right) \left( \bI_m + \sum_{i=1}^k \displaystyle \frac{\hat{\bW}_i }{s_j-\lambda_i} \right)^{-1} - \bH(s_j) \right)^2 .
\end{align}
This solution can, for example, be obtained from an optimization in least-squares sense. The weights for the SISO case are computed analogously. Here, the matrix-values weights and transfer function values reduce to scalars.

\subsubsection{Pole placement}
\label{sec:pole_placement}

The next step would be to take advantage of the degrees of freedom in the vector $\hat{\bB}$ from (\ref{ROM_SISO}), so that the ROM thus constructed has particular (stable) poles \cite{Ionescu2022,GarciaCanseco2015,Padoan2014}. These will be denoted with $\zeta_1, \zeta_2, \ldots, \zeta_k$. 
The following derivations assume a SISO model.
To enforce that this happens, we need to make sure that the matrix $\zeta_j \bI_k - \hat{\bA}$ loses rank for all $1 \leq j \leq k$.
In what follows, we show how to enforce this property in an elegant, straightforward way. Remember that the transfer function of the parameterized AF model is given by:
\begin{align}
	\hat{\bH}(s) &=  \frac{\sum_{i=1}^k \displaystyle \frac{\hat{w}_i h_i}{s-\lambda_i}}{1+\sum_{i=1}^k \displaystyle \frac{\hat{w}_i}{s-\lambda_i}} = \frac{N(s)}{D(s)}.
\end{align}
Now, let's say we would like this transfer function to have $k$ poles at the selected values $\zeta_j$'s. Clearly, the condition is $D(\zeta_j) = 0$ and hence we need to enforce:
\begin{align}
	&1+\sum_{i=1}^k \displaystyle \frac{\hat{w}_i}{\zeta_j-\lambda_i} = 0, \ \forall 1 \leq j \leq k\ \Leftrightarrow \
	 \cC_{\zeta,\lambda} \hat{\bB} = - \bone_k \ \Leftrightarrow \ \hat{\bB} = - \mathcal{C}_{\zeta,\lambda}^{-1} \bone_k,
\end{align}
where $\mathcal{C}_{\zeta,\lambda}$ is a Cauchy matrix defined by: $\left(\mathcal{C}_{\zeta,\lambda}\right)_{i,j} = \frac{1}{\zeta_i - \lambda_j}$.
Details on how to obtain the above expression by following the procedure in \cite{antoulas2007polplatzierung} are given in \Cref{sec:app_pp}. We note that placing poles is a difficult numerical problem which requires the inversion of a Cauchy matrix, which is highly ill-conditioned, by nature. 

Instead of doing this, we could solve $\mathcal{C}_{\zeta,\lambda} \hat{\bB} = - \bone_r$,
without inverting the Cauchy matrix explicitly, i.e., by solving a linear systems of equations.
\Cref{al:lfpp} summarizes this procedure in a data-driven context. The required underlying model is obtained from a set of transfer function evaluations by applying the Loewner framework. The method is illustrated for SISO systems, but can readily be extended to the MIMO case.

\begin{algorithm}[tb]
	\caption{Loewner framework with pole placement (\lfpp).}
	\label{al:lfpp}
	\begin{algorithmic}[1]
		\Require Transfer function samples $\left\{\bH\left(s_i\right)\right\}_{i=1}^N$, corresponding sampling points $\bXi = \left\{s_i\right\}_{i=1}^N$, locations for poles $\bm{\zeta} = \left\{\zeta_i\right\}_{i=1}^k$, interpolation points $\bm{\lambda} = \left\{\lambda_i\right\}_{i=1}^k$.
		\Ensure Surrogate model $\hat{\bH}(s) = \hat{\bC} (s \bI_r - \hat{\bA})^{-1} \hat{\bB}$.
		\State Compute $\Sig_D$ from $\left\{\bH\left(s_i\right)\right\}_{i=1}^N$ and $\left\{s_i\right\}_{i=1}^N$ using the Loewner framework (cf. \Cref{sec:LF}).
		\State $\hat{\bC} \gets \begin{bmatrix}
		\bH_D\left(\lambda_1\right) & \cdots & \bH_D\left(\lambda_k\right)
		\end{bmatrix}$.
		\State $\hat{\bB} \gets - \mathcal{C}_{\zeta,\lambda}\inv \bone_r$.
		\State $\hat{\bA} \gets \operatorname{diag}\left(\lambda_1, \ldots, \lambda_k\right) - \hat{\bB}\bone_k\trans $.
	\end{algorithmic}
\end{algorithm}

\subsection{Automatic choice of poles and interpolation points}
\label{sec:lfapp}

A reasonable choice of poles and interpolation points for \Cref{al:lfpp} is not always readily available, but the approximation of the surrogate is heavily influenced by this choice.
In the following, we show an extension to \Cref{al:lfpp} which computes a surrogate model \cref{eqn:Htrf_r} without requiring sets of poles and interpolation points as input parameters.
\Cref{al:lfapp} sketches the skeleton of such automatic algorithm. Similar to \Cref{al:lfpp} it employs the Loewner framework to obtain a realization of a surrogate interpolating the provided data.
Subsequently, a generalized eigendecomposition of the Loewner realization of the original data is computed to find suitable locations for poles.
From this is is possible to compute the dominance of all eigenvalues; for details, see, e.g. \cite{rommes2006efficient}. The algorithm now chooses the $k$ most dominant eigenvalues as poles to enforce in the surrogate. It should be noted that only eigenvalues with negative real parts should be considered, if the stability of the surrogate is important.
The required interpolation points can now be chosen similar to \Cref{al:ls_loewner} by computing a CUR decomposition and using the interpolation points associated to the rows or columns of the decomposition as interpolation points for the new surrogate.

The approximation of dominant poles of the underlying model from data is less robust if the transfer function samples are disturbed by noise. This leads to a reduced approximation quality. For a better performance if applied to noisy data, \Cref{al:lfapp} can be modified as follows:
To obtain poles which should be enforced, first choose manually the most prominent features in the transfer function, e.g. peaks, which should be approximated by the surrogate model.
Now choose the eigenvalues which imaginary parts are closest to the frequencies, where the chosen features of the transfer function are located.
The CUR decomposition also fails at extracting the most dominant rows and columns of the Loewner matrix if noisy data is assessed.
Therefore another heuristic is employed to choose the interpolation points:
Use the value $s_i$ which corresponds to the lowest amplitude of the transfer function between the locations of two enforced poles. This leads to reasonable approximations, especially for lightly damped systems. Other approaches include choosing simply the middle between the location of two poles or specifying an offset between pole and interpolation point location.

\begin{algorithm}[tb]
	\caption{Loewner framework with automatic pole placement (\lfapp).}
	\label{al:lfapp}
	\begin{algorithmic}[1]
		\Require Transfer function samples $\left\{\bH\left(s_i\right)\right\}_{i=1}^N$, corresponding sampling points $\bXi = \left\{s_i\right\}_{i=1}^N$.
		\Ensure Surrogate model $\hat{\bH}(s) = \hat{\bC} (s \bI_r - \hat{\bA})^{-1} \hat{\bB}$.
		\State Compute $\Sig_D$ from $\left\{\bH\left(s_i\right)\right\}_{i=1}^N$ and $\left\{s_i\right\}_{i=1}^N$ using the Loewner framework (cf. \Cref{sec:LF}).
		\State Compute the generalized eigenvalue decompositions $\bA \bX = \bE \bX \bm{\alpha}$ and $\bY\herm \bA = \bm{\alpha} \bY\herm \bE$ for the matrices of right and left eigenvectors $\bX, \bY$ and the matrix of eigenvalues $\bm{\alpha} = \operatorname{diag}\left(\alpha_1, \ldots, \alpha_n\right) $.
		\State Compute eigenvalue dominance $d_i = \frac{\lvert \bC \bY(:,i) \alpha_i \bX(:,i)\herm \bB \rvert}{\lvert \Re\left(\alpha_i\right) \rvert},\ i = 1, \ldots, n $ and sort $\bm{\alpha}$ accordingly
		\State Set $\bm{\zeta}$ to the $k$ most dominant eigenvalues.
		\State Compute CUR decomposition of $\IL$.
		\State Set $\bm{\lambda}$ to the $k$ right or left interpolation points corresponding to the CUR decomposition.
		\State Compute surrogate as in \Cref{al:lfpp}.
	\end{algorithmic}
\end{algorithm}

\section{Numerical results}
\label{sec:numerics}

In the following, we demonstrate the methods discussed in \Cref{sec:main} by applying them on three benchmark examples available from the MOR-Wiki\footnote{\url{http://modelreduction.org}}:
\begin{description}
	\item[ISS] This system models the structural response of the Russian Service Module of the International Space Station (ISS) \cite{morGugAB01}. The model has $n=270$ states, $m=3$ inputs, and $p=3$ outputs. The dataset used for the computations contains transfer function measurements at 400 logarithmically distributed points in the range $\left[10^{-1}, 10^2\right] \cdot \imath$.
	The model is also part of the SLICOT benchmark collection \cite{slicot_iss}.
	\item[Flexible aircraft] This system models lift and drag along the flexible wing of an aircraft. The system matrices are not available, we only have access to a dataset of 420 transfer functions samples at linearly distributed frequencies between 0.1 and \SI{42.0}{\Hz}. The original dataset has one input (the gust disturbance) and 92 outputs. For the following experiments, we choose the 91st output which corresponds to the first flexible mode~\cite{PoussotMATHMOD:2018}. The dataset is available from~\cite{MORWikiAircraft}.
	\item[Sound transmission] This system models the sound transmission through a system of two brass plates with an air enclosure between them. The transfer function measures the sound pressure in an adjacent acoustic cavity. The geometry is based on \cite{guy1981transmission}; the data---transfer function evaluations at 1000 linearly-distributed frequency values between 1 and \SI{1000}{\Hz}---is available from \cite{dataAum22}.
\end{description}

We note that no tangential interpolation (as described in \cite{MA07}) is applied for the MIMO model. Instead, the Loewner matrices are constructed in a block-wise manner.
The case of tangential interpolation, within the proposed approaches in this note, will be investigated in future works.

We enforce realness of all surrogate models (all matrices contain only real entries) by applying the transformation described in \cite{ALI17}. For this, all data must be available in complex conjugate pairs. The required transformation matrix is given by
\begin{equation}
	\bJ = \bI_\ell \otimes \left(\frac{1}{\sqrt{2}}\begin{bmatrix}
		\bI_m & \bI_m \\ -\imath\bI_m & \imath\bI_m
	\end{bmatrix}\right),
\end{equation}
with $\ell=\frac{k}{2}$ and the real-valued quantities are obtained from $\hat{\bA}^{(\Re)} = \bJ \hat{\bA} \bJ\herm$, $\hat{\bB}^{(\Re)} = \bJ \hat{\bB}$, and $\hat{\bC}^{(\Re)} = \hat{\bC} \bJ\herm$.

For some of the experiments we add artificial noise to the measurements, in order to obtain perturbed data. The modified measurements are given by
\begin{equation}
	\check{H}\left(s_i\right) = H\left(s_i\right)\left(1 + Z_i\right),\ i=1, \ldots, n,
\end{equation}
where $Z_i \in \IC$ is the $i$th sample drawn from a set of random numbers $Z \sim \cC\cN\left(\mu,\sigma^2\right)$ following a complex normal distribution with mean $\mu$ and standard deviation $\sigma^2$. Here, the real and imaginary parts of $Z$ are independent normally distributed variables \cite{drmac2022learning}.

We assess the approximation error of the surrogate models with an approximated $\cL_\infty$ norm, because many surrogates have unstable poles and hence, the $\cH_\infty$ can not be computed.
For a given reduced order $r$, the $\cL_\infty$ error in the considered frequency range $\omega \in [\omega_{\min},\omega_{\max}]$ is approximated by
\begin{equation}
	\varepsilon(r) = \frac{\max\limits_{\omega \in
			[\omega_{\min},\omega_{\max}]}
		\left\lVert H(\omega \imath) - \hat{H}_{r}(\omega \imath) \right\lVert_{2}}
	{\max\limits_{\omega \in [\omega_{\min},\omega_{\max}]}
		\left\lVert H(\omega \imath) \right\lVert_{2}}
	\approx \frac{\linfnorm{H - \hat{H}_{r}}}{\linfnorm{H}}.
\end{equation}
Note that strategies to post-process surrogates to obtain stable models have been studied in \cite{morGosPA21a}.

The numerical experiments have been conducted on a laptop equipped with an AMD Ryzen\texttrademark \\ ~7~PRO~5850U and
12\,GB RAM running Linux Mint 21 as operating system. All algorithms have been implemented and run with MATLAB R2021b Update~2 (9.11.0.1837725).

\begin{center}%
	\setlength{\fboxsep}{5pt}%
	\fbox{%
		\begin{minipage}{.92\linewidth} \small
			\textbf{Code and data availability}\newline
			The data that support the findings of this study are openly available in Zenodo at
			\vspace{.5\baselineskip}
			\begin{center}
				\href{https://doi.org/10.5281/zenodo.7490158}%
				{\texttt{doi:10.5281/zenodo.7490158}}
			\end{center}
			\vspace{.5\baselineskip}
			under the BSD-2-Clause license, authored by Quirin Aumann and
			Ion Victor Gosea.
	\end{minipage}}
\end{center}

\subsection{Case of exact measurement data}
In the following, we compare the performance of the new approach \AFAAALF{} to the following established strategies:
\begin{itemize}
	\item \LoewnerSVD: Truncate Loewner matrices populated with the complete dataset to order $r$ using an SVD~\cite{ALI17}.
	\item \AAACUR: Construct a purely interpolatory model of order $r$ using all data points chosen by the CUR decomposition, similar to \cite{morKarGA19a}.
	\item \AAA: Apply the strictly-proper variant of AAA \cite{morGosPA21a} to the complete dataset to compute a reduced-order model of size $r$.
\end{itemize}

\begin{figure}[tb]
	\centering
	\includegraphics{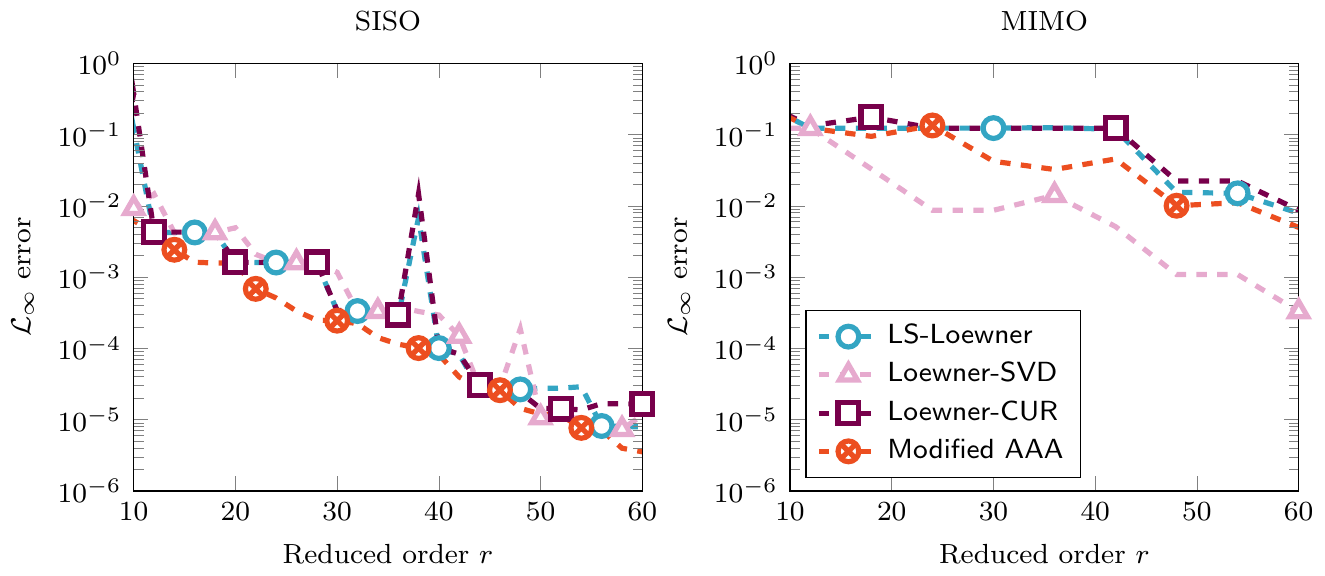}
	\caption{The approximated $\cL_\infty$ errors of reduced-order models of order $r$ computed from the ISS data. Left: SISO with first input and output, respectively. Right: MIMO with three inputs and three outputs $m=p=3$ (the number of interpolation points is $k=\frac{r}{m}$).}
	\label{fig:iss_hinf}
\end{figure}

We first consider the original MIMO ISS example and a SISO variant where we select the first input and output, respectively, from the MIMO system. To evaluate the overall performance of the different methods related to the size of a surrogate model, we compute the approximated $\cL_\infty$ errors for models with orders $6 \leq r \leq 60$.
The approximation error versus the dimension of the respective surrogate model is depicted in \Cref{fig:iss_hinf} for all four methods.

Since tangential interpolation was not employed here, the order of the MIMO surrogates rises by $m$ for each additional interpolation point, i.e., $r=km$.
This explains the lower accuracy of the MIMO surrogate. For the maximum reduced order $r=60$, $k=20$ interpolation points are considered. The errors of the SISO surrogates for $r=20$, i.e., $k=20$, is in a similar range as in the MIMO case.
The SISO surrogates reach similar levels of approximation for all employed methods. In the MIMO case, \LoewnerSVD{} performs best. This can be explained by the following observation: the other methods always consider the complete transfer function measurement $\bH(\lambda_i) \in \IC^{p \times m}$ per interpolation point, while \LoewnerSVD{} extracts only the $r$ most dominant singular vectors for projection, regardless of to which interpolation point they belong to. In turn, the other methods also consider probably less important parts of the data as long as one input/output combination of the respective sample is relevant for approximation.
It can also be noted that \AFAAALF{} and \AAACUR{} perform very similar. This was expected, as both methods rely on the same interpolation points.

All four methods are now employed to compute a surrogate model of size $r=108$ to approximate the transfer function of the flexible aircraft model. The size of the surrogate model is determined by truncating all singular values $\tau<\num{1e-6}$ of an underlying Loewner matrix.

\begin{figure}[tb]
	\centering
	\includegraphics{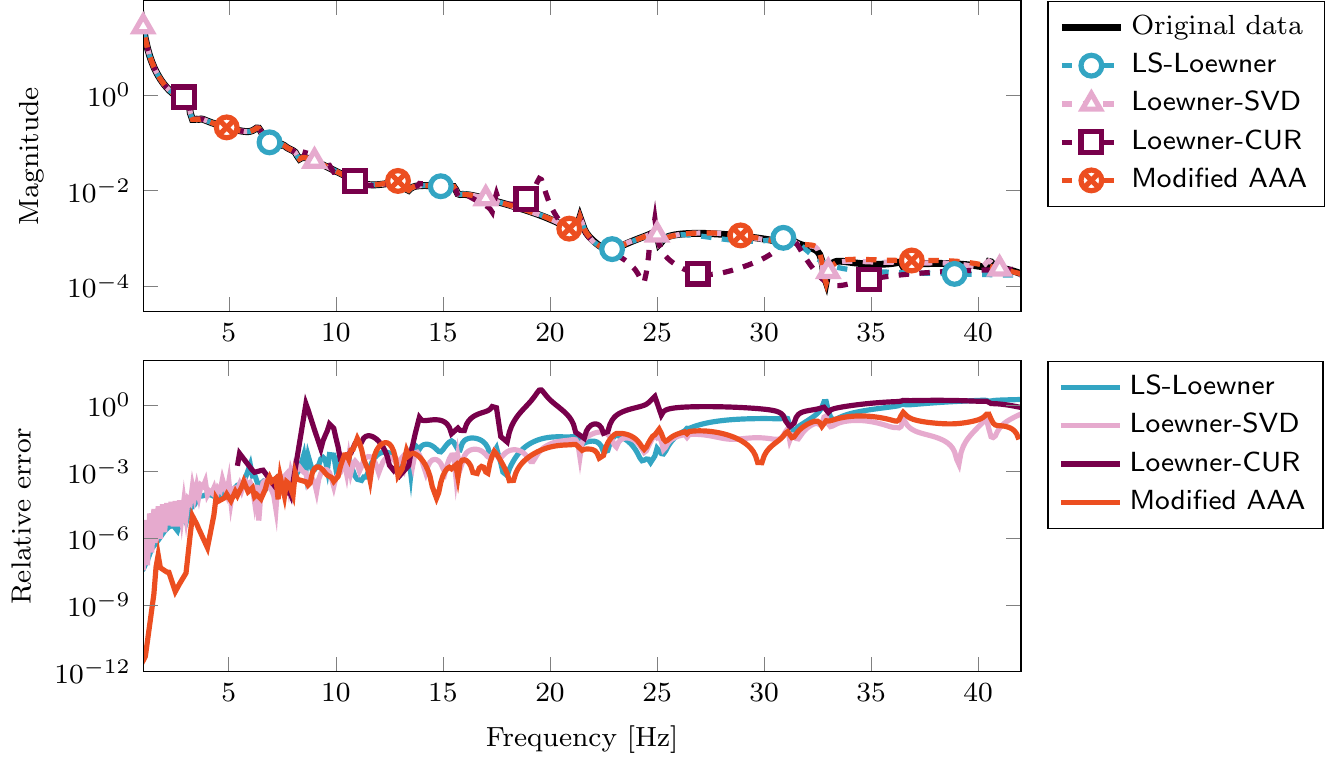}
	\caption{Transfer function (top) and relative pointwise errors (bottom) for reduced-order models of size $r=108$ for the aircraft model. The error is plotted only at frequencies which do not coincide to interpolation points of the respective method. }
	\label{fig:aircraft_tf}
\end{figure}

The transfer functions of all resulting models and their respective relative errors are given in \Cref{fig:aircraft_tf}.
Again, all methods succeed in computing a sufficiently accurate surrogate. However, the approximation quality of \AAACUR{} is noticeably worse than that of the other three methods.
Given that both \AAACUR{} and \AFAAALF{} use the same interpolation points, the weights computed from the least squares problem show a better performance compared to the partitioning approach used in \AAACUR.

\subsection{Perturbed measurement data}

\begin{figure}[tb]
	\centering
	\includegraphics{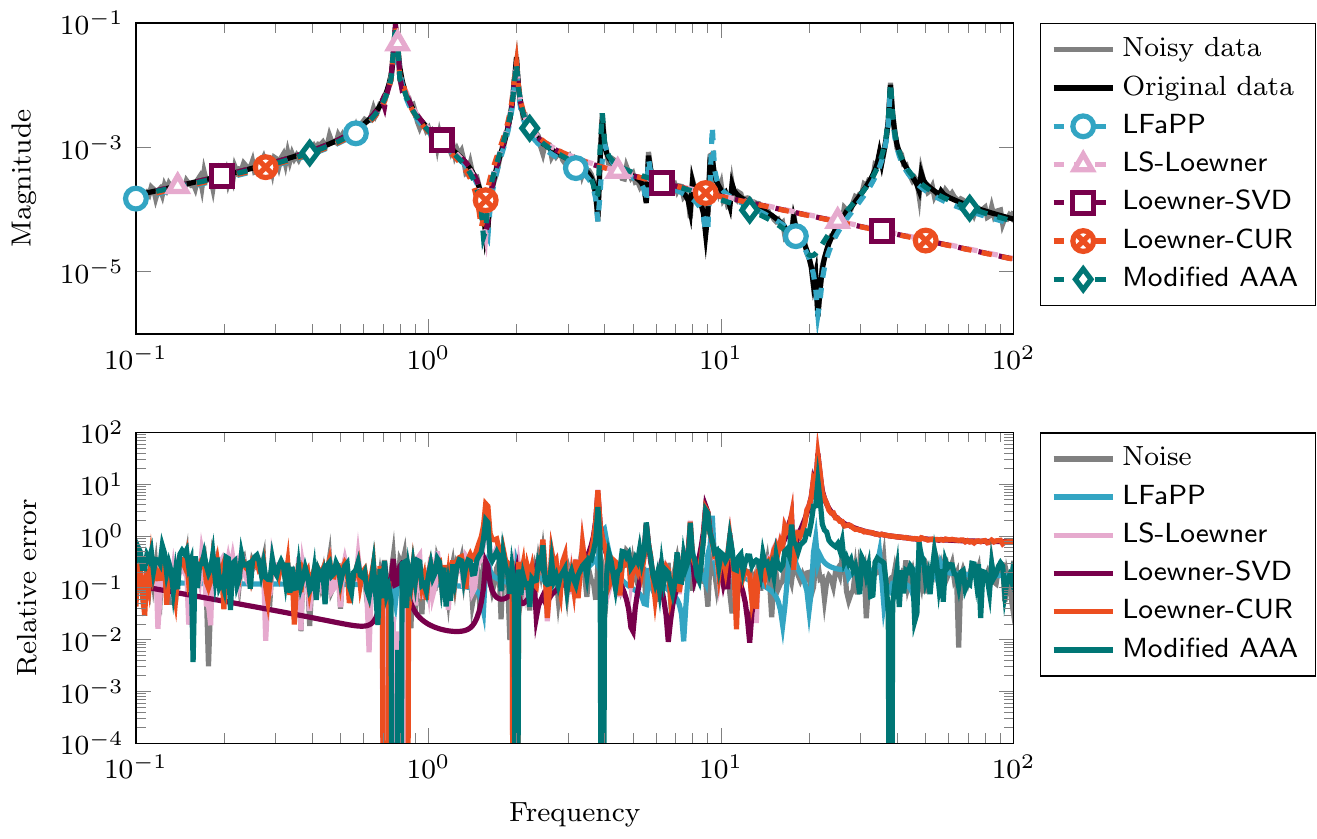}
	\caption{Transfer function of a surrogate with enforced poles compared to the noisy and original transfer function values. The transfer function of a model $r=12$ computed from \LoewnerSVD{} is given for reference.}
	\label{fig:iss_poles_tf}
\end{figure}

Analyzing measurement data perturbed by noise is a challenging task for interpolatory methods, such as the Loewner framework and the AAA algorithm (as pointed out in, e.g., \cite{GG20}).
In this experiment we investigate the effect of noise to the performance of the four methods described above and show, how enforcing poles and/or interpolation points can increase the approximation quality.
In the first experiment we consider transfer function data from the ISS model perturbed by noise with mean $\mu=0$ and standard deviation $\sigma^2=0.15$. We employ \lfapp{} and enforce poles at $\imath [.77, 2, 4, 5.6, 9.33, 37.9]$ near peaks of the transfer function. The resulting real-valued surrogate model has order $r=12$. The transfer functions of the surrogate model with enforced poles and reduced models computed from the same noisy data with \AFAAALF{}, \LoewnerSVD{}, \AAACUR, and \AAA{} are given in \Cref{fig:iss_poles_tf}. Enforcing the poles near peaks in the transfer function of the underlying data allows the surrogate to capture the behavior of the original data in a wider frequency range than applying  \AFAAALF{}, \LoewnerSVD{}, and \AAACUR. The choice of the locations, in which vicinity the poles should be chosen is, however, not automatized. \Cref{fig:iss_poles_tf} also shows the relative errors of all surrogate models referenced to the original data without noise.
While the enforced poles all have a negative real part, the models computed from the variants of the LF and AAA exhibit unstable eigenvalues. Thus, pole placement can be seen also as a means to enforce stability of the surrogate models. Alternatively, a post-processing step can be added to enforce stable models (for both LF and AAA methods), as performed in \cite{morGosPA21a}.

\begin{figure}[tb]
	\centering
	\includegraphics{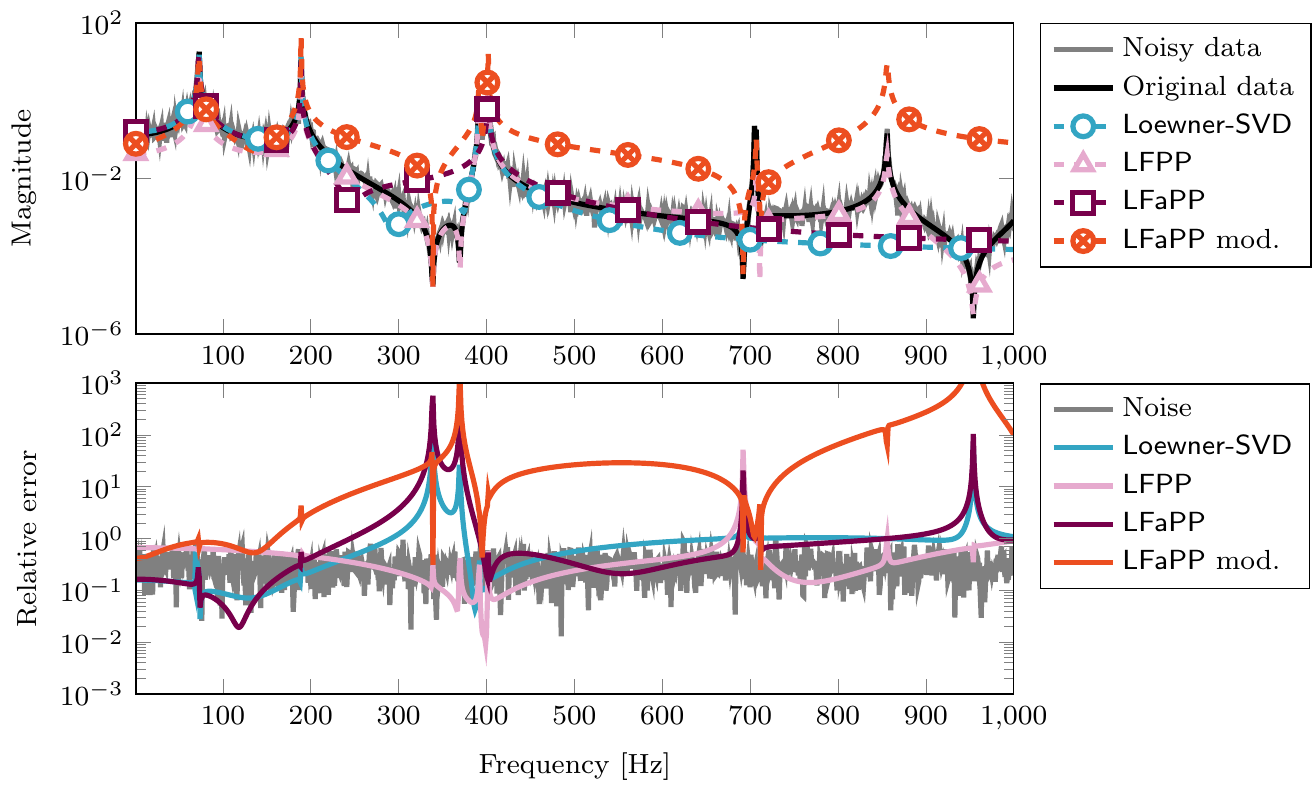}
	\caption{Transfer function (top) and relative pointwise errors (bottom) as well as the added noise for reduced-order models of size $r=12$ for the aircraft model. }
	\label{fig:transmission_cp_tf}
\end{figure}

We now evaluate the performance of the algorithms by applying them to heavily distorted transfer function measurements of the sound transmission problem.
Noise with a standard deviation of $\sigma^2=0.25$ is considered and
three algorithms are employed to compute surrogates: \LoewnerSVD, \lfpp{} (\Cref{al:lfpp}), and \lfapp{} (\Cref{al:lfapp}). We also test the modifications to \lfapp{} described in \Cref{sec:lfapp}. These results are denoted by \enquote{\lfapp{} mod.}.
For \lfpp{} we enforce poles at the eigenvalues of the underlying Loewner model which imaginary parts are near $2\pi\imath\left[72, 189, 392, 401, 706, 856\right]$. These locations correspond to characteristic peaks in the transfer function.
Further, we choose the interpolation points at $2\pi\imath\left[ 138, 339, 369, 569, 712, 954\right]$, which lie at the dips between the enforced poles.
\LoewnerSVD{} and \lfapp{} do not require input parameters in addition to the measured data.
\Cref{fig:transmission_cp_tf} shows the transfer function of the resulting surrogate models in comparison to the original and noisy underlying data.
It can be observed, that the automatic approaches \LoewnerSVD{} and \lfapp{} (mod.) cannot approximate the transfer function well after the first two peaks, i.e., for frequencies higher than \SI{200}{\Hz}, while \lfpp{} approximates the original data over the complete frequency range with decent accuracy.
The importance of reasonable interpolation points can be seen in the difference of \lfpp{} and \lfapp{}~mod., which have the same poles.
It should be noted that the surrogate model computed by \LoewnerSVD{} has two unstable poles while the other three surrogate models are stable.
It is, however, not always clear a priori how to choose the poles and interpolation points for \lfpp{} in order to achieve the best approximation quality possible. In this example, the noise level is too high for one of the automatic approaches to yield reasonable dominant interpolation points or poles.

\section{Conclusion and outlook}
\label{sec:conc}

In this contribution, we have proposed an extensive study of interpolation-based data-driven approaches for approximating the response of linear dynamical systems. All methods require input and output data, i.e., transfer function measurements, while direct access to the system operators or the states is not required. We showed different approaches how to achieve compact surrogate models approximating the input/output behavior of the original system and how to ensure various properties of the surrogate models, such as stability. Strategies how to work with noisy measurement data have also been addressed.

A natural extension of the framework described here is to apply the ideas of tangential interpolation as a means of modeling a MIMO system from data. Here, the tangential directions need to be incorporated in the parameterized one-sided realization.
Further topics include enforcing different structures of the original model in the surrogate model, e.g., second-order or delay structures.
It would also be interesting to study the possibility of placing certain stable poles while achieving interpolation in a least-squares sense.
Application cases for the proposed methodology could include damping optimization. Here, a family of parameterized interpolants could be used to find optimal positions for viscous dampers in a structural system.

\appendix
\section{Appendix}

\subsection{The Woodbury matrix identity}
\label{sec:app_wmi}

We can expand the right part of \cref{eq:realization_A}, such that:
\begin{equation}
	\hat{\bA} = \bLambda - \hat{\bB} \bR \Rightarrow s \bI_{km} - \hat{\bA}  = \underbrace{s \bI_{km} -\bLambda}_{\hat{\bM}} + \underbrace{\begin{bmatrix} \hat{\bW}_1 \\ \vdots \\ \hat{\bW}_k \end{bmatrix}}_{\hat{\bU}} \underbrace{\begin{bmatrix} \bI_m  \end{bmatrix}}_{\hat{\bT}} \underbrace{\begin{bmatrix} \bI_m & \cdots & \bI_m \end{bmatrix}}_{\hat{\bV}}.
\end{equation}
The \textbf{Woodbury matrix identity} is as follows:
$$
{\displaystyle \left(\hat{\bM}+ \hat{\bU}  \hat{\bT} \hat{\bV} \right)^{-1}=\hat{\bM}^{-1}-\hat{\bM}^{-1}  \hat{\bU} \left( \hat{\bT}^{-1}+ \hat{\bV} \hat{\bM}^{-1}  \hat{\bU} \right)^{-1}  \hat{\bV}  \hat{\bM}^{-1},}
$$
where $\hat{\bM}$, $\hat{\bU}$, $\hat{\bT}$ and $ \hat{\bV}$ are conformable matrices: $\hat{\bM}$ is $n \times n$, $\hat{\bT}$ is $k \times k$, $\hat{\bU}$ is $n \times k$, and $\hat{\bV}$ is $k\times n$. This can be derived using blockwise matrix inversion.
By denoting with $\bLambda_s  = s\bI_{km}-\bLambda$, then the first transfer function of the fitted model is written:
\begin{align}
	\begin{split}
		\hat{\bH}(s) &= \hat{\bC} \left( s \bI_{3m} - \hat{\bA} \right)^{-1} \hat{\bB} = \hat{\bC} \left( \bLambda_s +  \hat{\bU}  \hat{\bV} \right)^{-1} \hat{\bB} \\
		&=  \hat{\bC} \bLambda_s^{-1} \hat{\bB} - \hat{\bC} \bLambda_s^{-1}  \hat{\bU} \left(\bI_m+ \hat{\bV} \bLambda_s^{-1} \hat{\bU} \right)^{-1} \hat{\bV} \bLambda_s^{-1} \hat{\bB} \\
		&=  \hat{\bC} \bLambda_s^{-1} \hat{\bB} - \hat{\bC} \bLambda_s^{-1}  \hat{\bB} \left(\bI_m+ \bR \bLambda_s^{-1} \hat{\bB} \right)^{-1}\bR \bLambda_s^{-1} \hat{\bB} \\
		&= \hat{\bC} \bLambda_s^{-1} \hat{\bB} \left[ \bI_m - \left( \bI_m + \hat{\bX} \right)^{-1} \hat{\bX} \right] = \hat{\bC} \bLambda_s^{-1} \hat{\bB} \left( \bI_m + \hat{\bX} \right)^{-1},
	\end{split}
\end{align}
where $\hat{\bX} = \bR \bLambda_s^{-1} \hat{\bB}$. Hence, we arrive at \cref{eq:rom_woodbury} and the transfer function $\hat{\bH}(s)$ can be written as follows:
\begin{align}\tag{\ref{eq:rom_woodbury}}
	\hat{\bH}(s) &= \hat{\bC} \bLambda_s^{-1} \hat{\bB} \left( \bI_m + \bR \bLambda_s^{-1} \hat{\bB} \right)^{-1}.
\end{align}

\subsection{Pole placement as in \cite{antoulas2007polplatzierung}}
\label{sec:app_pp}

In order to enforce both prescribed poles and certain interpolation conditions in the ROM, we follow the derivations from \cite{antoulas2007polplatzierung}. It is to be noted that this approach is intrusive, i.e., requires access to the system's matrices. Hence, a descriptor model characterized in (generalized) state-space by the following equations
\begin{equation}\label{def_lin_sys_descriptor}
\Sig_{\textrm{Des}}: \begin{cases}
\bE \dot{\bx}(t) = \bA \bx(t) +\bB \bu(t), \ \
\by(t) = \bC \bx(t),
\end{cases}
\end{equation}
with corresponding transfer function ${\bH}_\textrm{Des}(s) = {\bC} (s {\bE} - {\bA})^{-1} {\bB}$ is considered to be given.
For the (right) interpolation points $\lambda_i,\ i=1, \ldots, k$ (where interpolation is imposed), and the desired poles to be placed, denoted with $\zeta_j$'s, the author in \cite{antoulas2007polplatzierung} starts by finding a row vector $\bC_{\zeta} \in \IC^{1 \times n}$ so that:
\begin{align}
\bC_{\zeta} \left[ \begin{matrix} (\lambda_1 \bE -\bA)^{-1} \bB \cdots (\lambda_k \bE -\bA)^{-1} \bB \end{matrix}  \right] = \bfz_{1 \times k}.
\end{align}
Then, the next step is to choose projection matrices $	\bW, \bV \in \IC^{n \times k}$ as
\begin{align}
\bW\herm = \begin{bmatrix}
\bC_{\zeta} (\zeta_1 \bE -\bA)^{-1} \\
\vdots \\
\tilde{\bC} (\zeta_k \bE -\bA)^{-1}
\end{bmatrix}, \ \
\bV = \left[ \begin{matrix} (\lambda_1 \bE -\bA)^{-1} \bB \cdots (\lambda_k \bE -\bA)^{-1} \bB \end{matrix}  \right].
\end{align}
As explained in \cite{antoulas2007polplatzierung}, the choice of $\bW\herm$ above is explained by imposing the required poles for the reduced model, while $\bV$ is chosen to match the interpolation conditions at the $\lambda_i$'s.
Moreover, using these notations, it follows that $\tilde{\bC} \bV = \bfz$. Next, put together the following matrices $\tilde{\bE} = \bW\herm \bE \bV, \ \ \tilde{\bA} = \bW\herm \bA \bV$. Then, it follows that $(s\tilde{\bE} - \tilde{\bA})$ loses rank when $s \in \{\zeta_1,\ldots,\zeta_r\}$. To show this, we simply write
\begin{align}
\bfe_j\trans (\zeta_j \tilde{\bE} - \tilde{\bA})  &= \bfe_j\trans \bW\herm (\zeta_j \bE - \bA) \bV=  \bC_{\zeta} (\zeta_j \bE - \bA)^{-1} (\zeta_j \bE - \bA) \bV =  \bC_{\zeta} \bV = \bfz.
\end{align}
Let $\bH_\zeta(s) =\bC_{\zeta} (s \bE -\bA)^{-1} \bB$ be a rational function in $s$ and we note that $\hat{\bE}$ and $\hat{\bA}$ are a special type of diagonally scaled Cauchy matrices, with the following exact definition:
\begin{align}\label{eq:defEAtilde}
\begin{split}
\tilde{\bE}_{i,j} &= -\frac{\bC_{\zeta} (\zeta_i \bE -\bA)^{-1} \bB - \bC_{\zeta} (\lambda_j \bE -\bA)^{-1} \bB}{\zeta_i - \lambda_j} = -\frac{\bH_\zeta(\zeta_i)}{\zeta_i - \lambda_j} \\
\tilde{\bA}_{i,j} &= -\frac{\zeta_i \bC_{\zeta} (\zeta_i \bE -\bA)^{-1} \bB - \lambda_j \bC_{\zeta} (\lambda_j \bE -\bA)^{-1} \bB}{\zeta_i - \lambda_j} = -\frac{\zeta_i \bH_\zeta(\zeta_i) }{\zeta_i - \lambda_j}
\end{split}
\end{align}
From the definition in \cref{eq:defEAtilde}, it follows that $\tilde{\bE} = -\bD_{\tilde{\bB}} \mathcal{C}_{\zeta,\lambda}$,
where $\bD_{\tilde{\bB}} = \text{diag}(\tilde{\bB})$ is a diagonal matrix. Similarly, it follows that $\tilde{\bA} = - \bZ \bD_{\tilde{\bB}} \mathcal{C}_{\zeta,\lambda}$, where $\bZ = \text{diag}(\zeta_1,\ldots,\zeta_k)$.

Next, we write the other projected quantities as
\begin{align}
\tilde{\bB} &= \bW\herm \bB = \begin{bmatrix}
\bH_\zeta(\zeta_1) & \cdots & \bH_\zeta(\zeta_k)
\end{bmatrix}\trans,  \ \
\tilde{\bC} = \bC \bV = \begin{bmatrix}
\bH(\lambda_1) & \cdots & \bH(\lambda_k)
\end{bmatrix} 
\end{align}
Hence, the reduced-order linear descriptor system $\Sigma_{\textrm{pp}}: (\tilde{\bE},\tilde{\bA},\tilde{\bB},\tilde{\bC})$ matches $k$ interpolation conditions and has the required poles.

Next, we show that this model can be written equivalently in the AF format. We first note that $\hat{\bC} = \tilde{\bC}$. For next step, provided that the matrix $\tilde{\bE}$ is non-singular, we remove it by incorporating it into the other matrices, as: $\breve{\bA} = \tilde{\bE}^{-1} \tilde{\bA}, \ \ \breve{\bB} = \tilde{\bE}^{-1} \tilde{\bB}, \ \
\breve{\bE} = \bI_k, \ \ \breve{\bC} = \tilde{\bC}$. We note that the two realizations of the interpolatory ROM, i.e., $(\hat{\bA}, \hat{\bB}, \hat{\bC})$ in \cref{ROM_SISO} and $(\breve{\bA}, \breve{\bB}, \breve{\bC})$ introduced above, are actually identical. The reason for this is that $\breve{\bC} = \hat{\bC}$ and the two ROMs match the same $k$ moments. Hence, it also follows that $\breve{\bB} = \hat{\bB}$.
Now, since $\breve{\bB} = \tilde{\bE}^{-1} \tilde{\bB}$ and $\tilde{\bE} = -\bD_{\tilde{\bB}} \mathcal{C}_{\zeta,\lambda}$, we can write that
\begin{align}
\hat{\bB} = - (\bD_{\tilde{\bB}} \mathcal{C}_{\zeta,\lambda})^{-1} \tilde{\bB} = - \mathcal{C}_{\zeta,\lambda}^{-1} \bD_{\tilde{\bB}}^{-1} \tilde{\bB} =- \mathcal{C}_{\zeta,\lambda}^{-1} \bone_k.
\end{align}
Hence, the above choice of vector $\hat{\bB}$ in \cref{ROM_SISO} imposes the required poles.


\addcontentsline{toc}{section}{References}
\bibliographystyle{plainurl}
\bibliography{authorref}

\end{document}